\documentclass{gtart}

\def\ifplaintex{\expandafter\ifx\csname documentclass\endcsname\relax}

\def\gtp{{\mathsurround=0pt\it $\cal G\mskip-2mu$eometry \&\ 
$\cal T\!\!$opology $\cal P\!$ublications}}  

\def\recd{{\small Received:\qua\receiveddate\ifx\reviseddate\relax
\else\qquad Revised:\qua\reviseddate\fi\par}} 

\def\lognumber#1{\def\thelognumber{#1}}
\def\volumenumber#1{\def\thevolumenumber{#1}}
\def\volumeyear#1{\def\thevolumeyear{#1}}
\def\papernumber#1{\def\thepapernumber{#1}}
\def\pagenumbers#1#2{\def\startpage{#1}\def\finishpage{#2}}
\def\published#1{\def\publishdate{#1}}

\def\received#1{\def\receiveddate{#1}}
\def\revised#1{\def\reviseddate{#1}}
\def\accepted#1{\def\accepteddate{#1}}

\def\asciiaddress#1{\def\theasciiaddress{#1}}

\long\def\asciiabstract#1{\long\def\theasciiabstract{#1}}

\let\\\par\let\thelognumber\relax\let\thevolumenumber\relax
\let\thepapernumber\relax\let\thevolumeyear\relax\let\startpage\relax
\let\finishpage\relax\let\publishdate\relax\let\receiveddate\relax
\let\reviseddate\relax\let\accepteddate\relax\let\theasciititle\relax
\let\theasciiauthors\relax\let\theasciiaddress\relax
\let\theasciiabstract\relax

\let\theasciiemail\relax

\ifplaintex
\font\logobig=cmssbx10 scaled 3836
\font\logomed=cmssbx10 scaled 2557
\else
\font\logobig=cmssbx10 scaled 4200
\font\logomed=cmssbx10 scaled 2800
\fi

\long\def\makeagttitle{   
\count0=\startpage
\agt\hfill      
\hbox to 45truept{\vbox to 0pt{\vglue -13truept{\logomed A\kern -.37em{\logobig 
T}\kern -.38em G}\vss}\hss}
\break
{\small Volume \thevolumenumber\ (\thevolumeyear)
\startpage--\finishpage\nl
Published: \publishdate}

\vglue .25truein

{\parskip=0pt\leftskip 0pt plus
1fil\def\\{\par\smallskip}{\Large\bf\thetitle}\par\medskip} \vglue
0.05truein

{\parskip=0pt\leftskip 0pt plus 1fil\def\\{\par}{\sc\theauthors}
\par\medskip}%
 
\vglue 0.03truein 

{\small\leftskip 25truept\rightskip 25truept{\bf Abstract}\stdspace\theabstract

{\bf AMS Classification}\stdspace\theprimaryclass
\ifx\thesecondaryclass\relax\else; \thesecondaryclass\fi\par
{\bf Keywords}\stdspace \thekeywords\par}\vglue 7truept

}   

\ifplaintex
\font\phead=cmsl9 scaled 950
\font\pnum=cmbx10 scaled 913
\font\pfoot=cmsl9 scaled 950
\headline{\vbox to 0pt{\vskip -4.5mm\line{\small\phead\ifnum
\count0=\startpage ISSN 1472-2739 (on-line) 1472-2747 (printed)
\hfill {\pnum\folio}\else\ifodd\count0\def\\{ }%
\ifx\theshorttitle\relax\thetitle\else\theshorttitle\fi\hfill{\pnum\folio}
\else\def\\{ and }{\pnum\folio}\hfill\ifx\theshortauthors\relax\theauthors
\else\theshortauthors\fi\fi\fi}\vss}}
\footline{\vbox to 0pt{\vglue 0mm\line{\small\pfoot\ifnum\count0=\startpage
\copyright\ \gtp\hfill\else
\agt, Volume \thevolumenumber\ (\thevolumeyear)\hfill\fi}\vss}}
\else
\font=cmsl9 scaled 1050
\font\lnum=cmbx10 
\font=cmsl9 scaled 1050
\makeatletter
\def\@oddhead{{\small\lhead\ifnum\count0=\startpage ISSN 1472-2739 
(on-line) 1472-2747 (printed)\hfill {\lnum\number\count0}\else\ifodd\count0
\def\\{ }\ifx\theshorttitle\relax \thetitle \else\theshorttitle\fi\hfill
{\lnum\number\count0}\else\def\\{ and }{\lnum\number\count0}
\hfill\ifx\theshortauthors\relax 
\theauthors\else\theshortauthors\fi\fi\fi}}\def\@evenhead{\@oddhead}
\def\@oddfoot{\small\lfoot\ifnum\count0=\startpage\copyright\ \gtp\hfill\else
\agt, Volume \thevolumenumber\ (\thevolumeyear)\hfill\fi}
\def\@evenfoot{\@oddfoot}
\makeatother
\fi
\let\maketitlepage\makeagttitle
\let\makeshorttitle\maketitlepage
\let\maketitle\maketitlepage

\newwrite\gtoutfile
\long\gdef\makeheadfile{  
{\def\\{, }\def\s{ }
\immediate\openout\gtoutfile head.xxx
\immediate\write\gtoutfile{To: math@arxiv.org}
\immediate\write\gtoutfile{Subject: put OR rep NNNNN:ppppp}
\immediate\write\gtoutfile{--text follows this line--}
\immediate\write\gtoutfile{Proxy-for: \ifx\theasciiauthors\relax
\theauthors\else\theasciiauthors\fi\s<\ifx\theasciiemail\relax\theemail\else\theasciiemail\fi>}
\immediate\write\gtoutfile{\noexpand\\}
\immediate\write\gtoutfile{Authors: \ifx\theasciiauthors\relax
\theauthors\else\theasciiauthors\fi}
{\def\\{ }\immediate\write\gtoutfile{Title: \ifx\theasciititle\relax
\thetitle\else\theasciititle\fi}}
\immediate\write\gtoutfile{Subj-class: GT or SG, GR etc}
\immediate\write\gtoutfile{MSC-class: \theprimaryclass\ifx\thesecondaryclass\relax\else, \thesecondaryclass\fi}
\immediate\write\gtoutfile{Journal-ref: Algebr. Geom. Topol. \thevolumenumber\s
(\thevolumeyear) \startpage-\finishpage}
\immediate\write\gtoutfile{Comments: Published by Algebraic and
Geometric Topology at}
\immediate\write\gtoutfile{\s\s\s  http://www.maths.warwick.ac.uk/agt/AGTVol\thevolumenumber/agt-\thevolumenumber-\thepapernumber.abs.html}
\immediate\write\gtoutfile{\noexpand\\}
\immediate\write\gtoutfile{}
\ifx\theasciiabstract\relax
\immediate\write\gtoutfile{\theabstract}\else
\immediate\write\gtoutfile{\theasciiabstract}\fi
\immediate\write\gtoutfile{}
\immediate\write\gtoutfile{\noexpand\\}
\immediate\write\gtoutfile{}
\immediate\closeout\gtoutfile}}  

\def\maketitlepage{\makeagttitle\makeheadfile}
\let\makeshorttitle\maketitlepage
\let\maketitle\maketitlepage

\lognumber{32}
\volumenumber{3}
\volumeyear{2003}
\papernumber{32}
\published{4 October 2003}
\pagenumbers{969}{992}
\received{29 January 2003}
\revised{14 August 2003}
\accepted{21 September 2003}

\usepackage{amssymb, amsmath}
\usepackage{latexsym}
\usepackage{epsf}

\newtheorem{thm}{Theorem}[section] 
\newtheorem{pro}[thm]{Proposition}
\newtheorem{lem}[thm]{Lemma}

\newtheorem{con}[thm]{Conjecture}

\theoremstyle{definition}
\newtheorem{dfn}{Definition}[section]

\def\1{{\rm1\mathchoice{\kern-0.25em}{\kern-0.25em}
        {\kern-0.2em}{\kern-0.2em}I}}

\newcommand{\lmn}[1]{\vadjust{\setbox1=\vtop{\hsize 25mm
\parindent=0pt\baselineskip=9pt
\rightskip=4mm plus 4mm#1}
\hbox{\kern-26mm\smash{\raise .5ex\box1}}}}

\newcommand{\psdiag}[3]{\hspace{1mm}\raisebox{-#1mm}{\epsfysize#2mm
\epsffile{#3.eps}}\hspace{1mm}}

\newcommand{\nc}{\newcommand}
\def\be#1\ee{\begin{equation}#1\end{equation}}
\nc{\bc}{\begin{center}}
\nc{\ec}{\end{center}}
\nc{\N}{{\mathsf N}}
\nc{\K}{{\mathsf K}}
\nc{\fk}{\mathbf{k}}
\theoremstyle{remark}
\newtheorem{rem}[thm]{Remark}

\def\v8{\vskip 8pt}
\def\a{\alpha}
\def\la{\langle}
\def\ra{\rangle}
\def\l{\lambda}
\def\n{\nu}
\def\m{\mu}

\begin{document}

\title[Geometric construction of spinors]{Geometric construction
of spinors\\in orthogonal modular categories}
\author{Anna Beliakova}

\address{Mathematisches Institut, Universit\"at Basel\\Rheinsprung 21,
CH-4051 Basel, Switzerland}
\asciiaddress{Mathematisches Institut, Universitaet Basel\\Rheinsprung 21,
CH-4051 Basel, Switzerland}

\email{Anna.Beliakova@unibas.ch}

\begin{abstract}
A geometric construction of ${\mathbb Z}_2$--graded 
odd and even orthogonal modular categories is given.
Their 0--graded parts coincide with categories previously
obtained by  Blanchet and the author  from the category of tangles
modulo the Kauffman skein relations. Quantum 
dimensions and twist coefficients of
1--graded simple objects (spinors)
are calculated. We show that invariants coming from our odd and even
orthogonal modular categories
admit spin and 
${\mathbb Z}_2$--cohomological refinements, respectively.
The relation with the quantum group approach is
discussed. \end{abstract}

\asciiabstract{A geometric construction of Z_2-graded odd and even 
orthogonal modular categories is given.  Their 0-graded parts coincide
with categories previously obtained by Blanchet and the author from
the category of tangles modulo the Kauffman skein relations. Quantum
dimensions and twist coefficients of 1-graded simple objects (spinors)
are calculated. We show that invariants coming from our odd and even
orthogonal modular categories admit spin and Z_2-cohomological
refinements, respectively.  The relation with the quantum group
approach is discussed.}

\maketitle

\section*{Introduction}
In 1993, Lickorish  gave a simple geometric construction
of 3--manifold invariants  based on the Kauffman brackets.
The same invariants were obtained earlier by Reshetikhin
and Turaev from the representation category of
the quantum group  ${ U}_q({\mathfrak s\mathfrak l}_2)$.
The method of Lickorish was so much easier than the quantum
group theoretical one that it inspired many researchers
to work on its generalizations.

Recall that a quantum group ${U}_q({\mathfrak g})$  
for {\it any} semi--simple Lie algebra ${\mathfrak g}$ 
and some root of unity $q$
provides 3--manifold invariants. In many cases
a representation category of the  quantum group  is modular
(or modularizable), i.e.\
 a  functor  from 
 the category of 3--cobordisms to the representation category
 can be constructed (see \cite{Le}).
This functor is called  a Topological Quantum Field Theory (TQFT). 

In order to get a geometric construction of 3--manifold invariants or TQFT's
coming from  quantum groups of type $A$ ($\mathfrak g= \mathfrak s
\mathfrak l_m$),  a replacement of
the Kauffman bracket in Lickorish's approach
by the HOMFLY polynomial is needed. This was successfully
done in papers of Yokota, Aiston--Morton and Blanchet.

In 
\cite{BB}, Blanchet and the author
 constructed (pre--)modular categories from the category of
tangles modulo the Kauffman skein relations. We recovered
the  invariants of symplectic
 quantum groups ($\mathfrak g= \mathfrak s
\mathfrak p_m$ type $C$), but only
 ``half'' of the invariants for orthogonal groups
($\mathfrak g= \mathfrak s
\mathfrak o_m$ types $B$ and $D$).
Our approach did not provide objects corresponding to spin representations.

In this article  we give a geometric construction of 
two series of orthogonal modular 
categories which include spinors.
We consider the  category  of framed tangles  
where colors from the set $\{1,2\}$
are attached to
 lines. We add the relations given 
by the kernel of the $({\mathfrak s\mathfrak o}_m, V, S)$ 
weight system pulled back by the Vassiliev--Kontsevich invariant.
The standard representation $V$ and the spinor representation $S$ 
are used for $1$--colored and  $2$--colored lines, respectively.
The resulting category admits a natural ${\mathbb Z}_2$--grading. 
The 0--graded part has the same set of simple objects as the
 category studied in  \cite{BB}. 

We consider two series of parameter specializations
 for which this 0--graded part is pre--modular.
Then  we  give a recursive construction of idempotents
for the 1--graded parts of these pre--modular categories. 
The key point  is the observation
that encircling any  1--graded object
with  a line colored with a special  0--graded object, given by 
a rectangular Young diagram,  yields a projector.
We calculate  quantum 
dimensions and twist coefficients of
1--graded simple objects (spinors).
We show that invariants coming from our odd and even
orthogonal modular categories
admit spin and 
${\mathbb Z}_2$--cohomological refinements, respectively.

This paper  is organized as follows. In the first 
section we define the category we will work in.
In the second and third sections we construct odd and even orthogonal 
modular categories. Relations with quantum groups are 
discussed in the last section.

The author wishes to thank Christian Blanchet for many stimulating
discussions.

\section{Basic category }
In this section we define the category which will 
be studied subsequently. We
also  analyze its 0--graded part.

\subsection{Category of two colored tangles}
Let us fix an oriented 3--dimensional Euclidean space ${\mathbb R}^3$ 
with coordinates $(x,y,t)$.
\begin{dfn}
A two colored tangle $T$ is 
a 1--dimensional compact smooth sub--manifold of ${\mathbb R}^3$  equipped 
with a 
normal vector field and lying between two horizontal planes $\{t=a\}$, 
$\{t=b\}$, $a>b$, called the top and the bottom planes.
The boundary of $T$ lies on two lines $\{t=a, y=0\}$ and $\{t=b, y=0\}$.
The normal vector field has coordinates $(0,1,0)$ in boundary points.
The  map $c$ from the set of connected components of $T$ to 
the set of colors $\{1,2\}$ is given. 

Two colored tangles
$T$ and $T'$ are equivalent if there is an isotopy
sending $T$ to $T'$ which respects horizontal planes and colorings.
\end{dfn}

Connected components of $T$ 
will be called lines. We represent $T$
by drawing its   generic position diagram in
blackboard framing. Lines of the second color are drawn  bold.
An example is given in Figure 1.

\begin{figure}
$$\psdiag{9}{27}{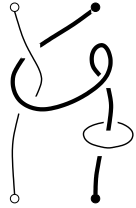}$$
\caption{Diagram of a two colored tangle}
\end{figure}

An intersection of a  two colored tangle with  the
top and the bottom planes defines a word in the alphabet $\{\circ, \bullet\}$,
where $\circ$ and $\bullet$ denote the points of the first and the second
color, respectively. For two such words $u$ and $v$, let $({\cal T},u,v)$
be the set of two colored tangles  whose intersection
with the top and the bottom planes are given by $u$ and $v$, respectively.

\begin{dfn}
Let ${\cal T}$ be the monoidal category whose objects are 
words in the alphabet $\{\circ, \bullet\}$.
For $u,v\in Ob({\cal T})$,  the set of  morphisms
$Hom(u,v)$ from $u$ to $v$
is  given by $({\cal T},u,v)$. The composition of $({\cal T},u,v)$ 
with $({\cal T},v,w)$ is defined
by gluing of horizontal planes identifying 
points corresponding to $v$. Moreover, $u\otimes v:= uv$.
\end{dfn}

\begin{dfn}
Let $f$ be a field. Let ${\cal T}_f$ be a linearization 
of $\cal T$, where formal $f$--linear combinations 
of tangles are allowed as morphisms.
The composition and tensor product are bilinear.
\end{dfn}

\subsection{Kontsevich integral}
In \cite{LM},  the category of q--tangles was considered.
The objects of this category are non--associative words in 
the alphabet  $\{+,-\}$. The morphisms are framed oriented tangles.
It was shown that  the universal Vassiliev--Kontsevich invariant
extends to a functor  from this category 
to the category of chord diagrams. An analogous construction 
applies to  colored q--tangles.

Let us orient all lines of two colored tangles from the top to the bottom.
Let us  map  a word  $u\in Ob({\cal T})$  with $n$ letters 
to the non--associative word of length $n$
 in  the alphabet $\{+,$  {\bf\bf +}$\}$ 
beginning with $n$
 left brackets,
 e.g.,  $\circ\bullet\circ$ maps to
$((+${\bf\bf +}$)+)$.
This defines a functor from  ${\cal T}$ into the category of
two colored q--tangles.  Now the universal Vassiliev--Kontsevich invariant
constructed in  \cite{LM} defines  a  functor  from this category to
the category $\cal A$ of chord diagrams with two colored support. 
We denote by $Z: {\cal T}\to {\cal A}$ the composition.

Let us consider the Lie algebra $\mathfrak g=\mathfrak s\mathfrak o_m$.
Let $V$ be the standard and $S$ 
the spin representation of $\mathfrak g$.  Let $t\in {\mathfrak g}^\ast \otimes {\mathfrak g}^\ast$
be its Killing form.

\begin{thm}
There exists a unique ${\mathbb C}$--linear monoidal 
 functor ${\cal F}_{\mathfrak g, V,S}$ (called weight system)
from $\cal A$ to the category  
 $Mod_{\mathfrak g}$ of the representations of $\mathfrak g$ 
such that it is uniquely characterized by its values on the following
elementary morphisms.

$$\psdiag{16}{48}{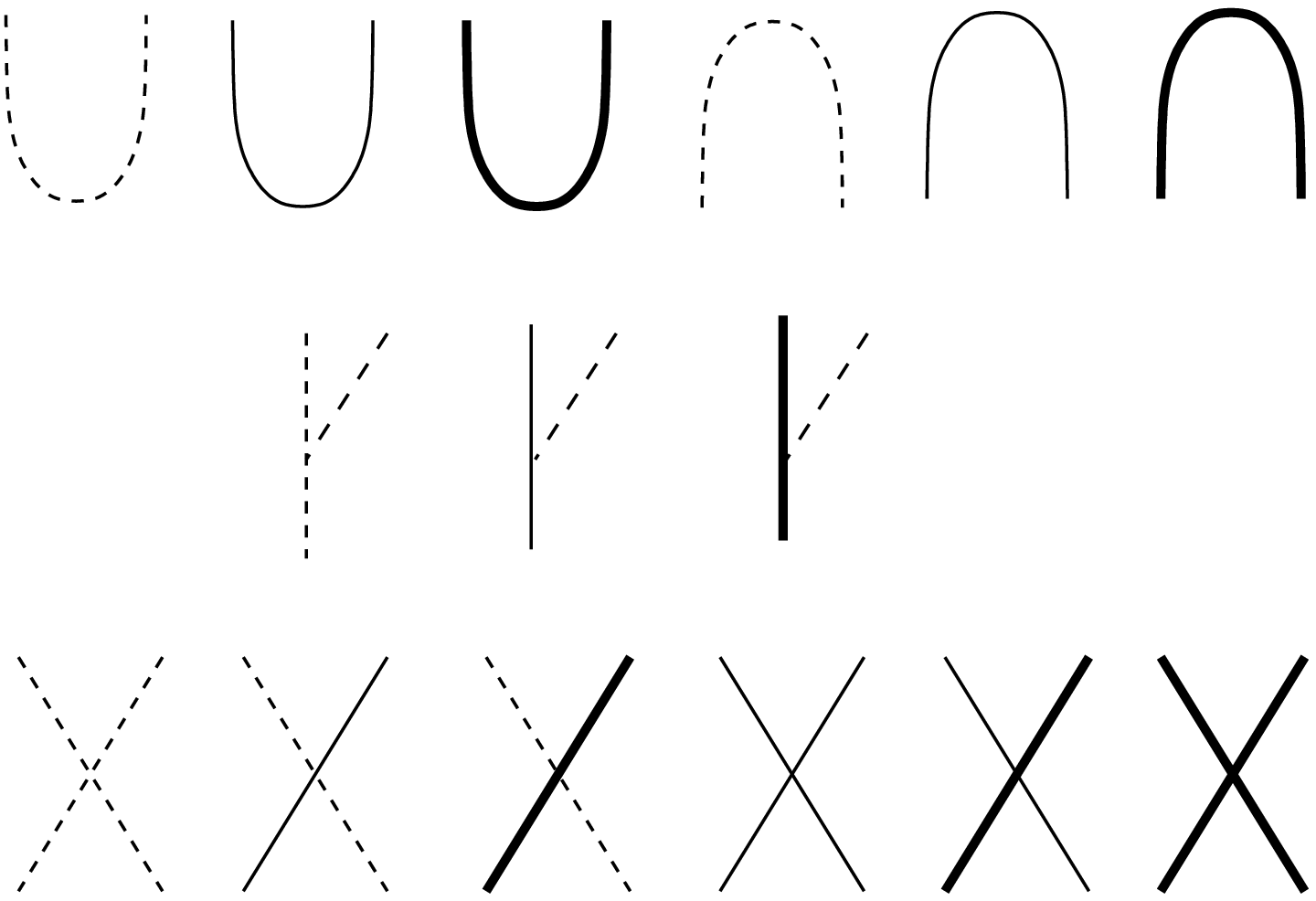}$$

\v8
The first three diagrams correspond 
to the  morphisms $\mathfrak g\otimes \mathfrak g \to {\mathbb C}$, 
$ V\otimes V\to {\mathbb C} $ and  $ 
S\otimes S\to {\mathbb C}$ in $Mod_{\mathfrak g}$  given by the Killing form. 
The next three are their transposes.
The first morphism in the second row is given by the Lie bracket.
The next two
correspond to the $\mathfrak g$--action on $V$ and $S$. The third row
describes  flips $x\otimes y \mapsto y\otimes x$ in 
$\mathfrak g\otimes \mathfrak g$, $\mathfrak g \otimes V$, 
$\mathfrak g \otimes S$,
$V\otimes V$, $V\otimes S$ and $S\otimes S$, respectively. 

\end{thm}

The proof of Bar--Natan  \cite{BN} can be adapted. An essential 
point is that the 
invariant tensors  $\varrho(t)$ for any representation $\varrho$ satisfy
the classical Yang--Baxter equation, which 
corresponds to the $4$--term relation in $\cal A$.

\medskip
{\bf Remark}\qua For $\mathfrak g =\mathfrak s\mathfrak o_{2n}$, 
the construction can be modified
by orienting  2--colored lines and by using the spin 
representations $S_\pm$ in the weight
system, according to the orientation.

\subsection{Category ${\cal T}_q(\mathfrak s\mathfrak o_m)$}
The central object of our study is the  category 
${\cal T}_h(\mathfrak s\mathfrak o_m)$
defined as ${\cal T}_{\mathbb C}$ modulo the relations given by the kernel of 
${\cal F}_{\mathfrak s\mathfrak o_m, V, S} ( Z({\cal T}_{\mathbb C}))$. 
The relations are defined {\it a priori} over $\mathbb C[[h]]$, where 
$h$ is the formal parameter of the Kontsevich integral. 
An explicite description of the relations is not known except 
if  we restrict to 1--colored framed tangles 
 or  to 2--colored ones and use  ${\cal F}_{\mathfrak s\mathfrak o_7, S}$
weight system. The first case was considered in \cite{LM1}
and the relations are just the Kauffman skein relations.
The second case was studied in \cite{P}, where a set 
of relations sufficient to calculate link invariants is given.

The proof of  Le and Murakami in \cite{LM} can be used to show
 that  link invariants provided by ${\cal T}_h (\mathfrak s\mathfrak o_m)$
and the quantum group $U_q(\mathfrak s\mathfrak o_m)$ coincide for odd $m$
if  $q=\exp h$ 
 and for even $m$ if $q=\exp 2h$.
The invariant associated by ${U}_q (\mathfrak s\mathfrak o_m)$
with a  colored  link is defined  over the ring $R={\mathbb Q}[q^{\pm
\frac{ 1}{2D}}]$, where $D=2$ if $m$ is odd and $D=4$ for even 
$m$.\footnote{The integrality result of T. Le shows that even a 
smaller ring can be considered.}
This allows  to define  ${\cal T}_h (\mathfrak s\mathfrak o_m)$ 
over $R$ and to 
 use the notation  ${\cal T}_q (\mathfrak s\mathfrak o_m)$
for  ${\cal T}_h (\mathfrak s\mathfrak o_m)$,  where
  $q=\exp h$ if $m=2n+1$ and  $q=\exp 2h$ if $m=2n$.

Let us define a ${\mathbb Z}_2$--grading in ${\cal T}_q (\mathfrak 
s\mathfrak o_m)$
as follows. A grading of $u\in Ob({\cal T}_q (\mathfrak s\mathfrak o_m))$
 is given by the number of symbols 
$\bullet$  in $u$ modulo 2.  All morphisms in the category are 0--graded.

\subsection{0--graded idempotents }

Let $u\in Ob({\cal T}_q (\mathfrak s\mathfrak o_m))$. A
nonzero morphism  $T\in
End_{{\cal T}_q (\mathfrak s\mathfrak o_m)}(u)$ is called a
minimal idempotent if
$T^2=T$ and for any $X\in End_{{\cal T}_q (\mathfrak s\mathfrak o_m)}(u)$ 
there exists a constant $c\in R$, such that  $TXT= c T$.
A standard procedure called idempotent completion allows
to add  idempotents as objects into the category. 
Objects given by minimal idempotents are called simple.
The idempotent completion of ${\cal T}_q (\mathfrak s\mathfrak o_m)$ is 
denoted by the same symbol. 
Let us equip ${\cal T}_q (\mathfrak s\mathfrak o_m)$ with a direct sum of objects in a formal way.

In \cite{BB1}, we gave a geometric  construction of 
minimal idempotents in the 
category of (framed non--oriented) tangles modulo the Kauffman skein relations.
\begin{figure}[h]
$$
\psdiag{3}{9}{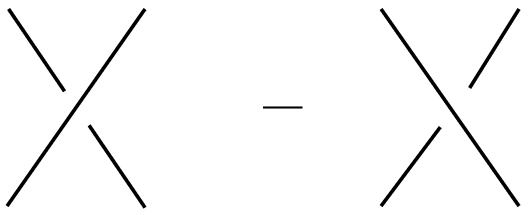} = \;(s-s^{-1})\;
\left(\;\,\psdiag{3}{9}{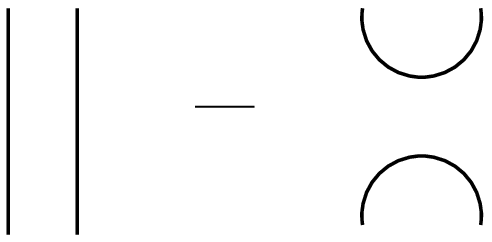}\;\,\right) 
$$
$$
\psdiag{3}{9}{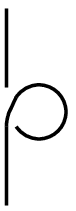}\;=\;\a\;\; \psdiag{3}{9}{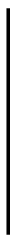}\;,\;\;\;\;\;
\psdiag{3}{9}{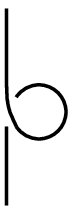}\;=\;\a^{-1} \;\;\psdiag{3}{9}{1_bel}$$

$$L\;\amalg \;\bigcirc\; = \;\left(\frac{\a-\a^{-1}}{s-s^{-1}}+1
\right)\;\, L  $$
\caption{\label{kaufrel} Kauffman relations}
\end{figure}

The idempotents were numbered 
 by integer partitions or Young diagrams
 $\l=(\l_1\geq\l_2\geq ...\geq\l_k)$.
Their twist coefficients and quantum dimensions were calculated.

\begin{lem}
After the substitution $\alpha=s^{m-1}$, $s=\exp h$,
the idempotents constructed in \cite{BB1} give the whole set
of minimal idempotents of 
the 0--graded part of ${\cal T}_h (\mathfrak s\mathfrak o_m)$.
\end{lem}
\begin{proof}
Let us first consider 1--colored tangles.
The set of relations  in  ${\cal T}_h (\mathfrak s\mathfrak o_m)$
for 1--colored tangles coincides with
 the Kauffman skein relations, where $\a=s^{m-1}$, $s=\exp h$ 
 (see \cite{LM1} for the proof, -- 
attention, --   Le and Murakami use a different normalization  for
the trivial knot). Minimal idempotents for  tangles modulo the 
Kauffman skein relations   are constructed in \cite{BB1}.

From the representation theory of the  classical orthogonal Lie algebras
(see e.g.\ \cite{FH} p. 291--296) we know that $S\otimes S$ decomposes
into a direct sum of simple objects numbered by integer partitions.
This implies that  the addition of an even number of 2--colored lines
to 1--colored tangles does not create new minimal idempotents.
\end{proof}

In \cite{BB}, we found seven series of  specializations of parameters
$\a$ and $s$, such that the category of tangles modulo the Kauffman
skein relations with these specia\-li\-zations becomes pre--modular
(after idempotent completion and quotienting by negligible morphisms).
In all  cases, $s$ is a root of unity and $\a$ is a power of $s$.
The specializations $\a=s^{2n}$, $s^{4n+4k}=1$ and
$\a=s^{2n-1}$, $s^{4n+4k-4}=1$ lead to odd and even orthogonal
categories  $B^{n,{-k}}$ and $D^{n,k}$, respectively.

In the remainder of the paper we will complete
 $B^{n,{-k}}$ and $D^{n,k}$ with
1--graded simple objects.
The  odd and even orthogonal  cases will be treated separately.

\section{Odd orthogonal modular categories}
This section is devoted to the construction of odd orthogonal
modular categories. We show that these
categories lead to spin TQFT's and calculate spin Verlinde formulas.

\subsection{0--graded objects}
Let us use the standard notation $B_n$ for
 $\mathfrak s\mathfrak o_{2n+1}$. We fix a 
primitive $(4n+4k)$th root of unity $q$ and choose $v$ with $v^2=q$.
In this specialization,
the set of 0--graded simple objects of
${\cal T}_q (B_n)$ (modulo negligible morphisms)
is given by Young diagrams (or integer partitions) from
the set  
$$\tilde \Gamma =\{\l:\l_1+\l_2\leq 2k+1,\l^\vee_1+\l^\vee_2\leq 2n+1\}$$
(see \cite{BB} p. 487 for the proof, put $s=q$). 
 Here $\l^\vee_i$ denotes the 
number of cells in the $i$th column of $\l$.

The set $\tilde\Gamma$ admits an algebra  structure
with the multiplication given by the tensor product and the addition 
given by the direct sum. 
The empty partition
is the one in this algebra and it will be denoted by $1$. 
(Please not confuse with the partition $V=(1)$ corresponding
to the object $\circ$.) 

The set $\tilde \Gamma$ has two invertible objects of order two.
The object $J=(2k+1)$, given by the  one row Young diagram with
$2k+1$ cells in it,  and the object $1^{2n+1}$ given by the
one column Young diagram with $2n+1$ cells in it.
They are  0--transparent, i.e.\ they have  trivial braiding 
with any other 0--graded object.
The tensor square of each of them is the trivial object.
The twist coefficient of $J$ is minus one and of $1^{2n+1}$ is one. 
Let us consider the following set
$$\Gamma_0=\{\l:\l_1+\l_2\leq 2k+1,\l^\vee_1\leq n\} .$$
Any $\l\in \tilde\Gamma$ is either contained in $\Gamma_0$
or is isomorphic to $1^{2n+1}\otimes  \m$ with $\m\in \Gamma_0$,
where $1^{2n+1}\otimes \m $ and $\m$ are both simple
with the same quantum dimensions,
braiding and twist  coefficients. 
There exists a standard procedure
called modularization (or modular extension),
which allows to path to a new  category, where 
$1^{2n+1}\otimes \m$ and $\m$ are identified. 
We will denote by  ${\cal T}_q (\tilde B_n)$ this new category.
Its set of simple 0--graded objects is $\Gamma_0$.
The existence criterion for
such  modularization functors was developed by
Brugui\`eres. 
In \cite{BB}, a  geometric construction of these functors is given.

\vspace*{0.3cm}

\subsection{Recursive construction of 1--graded idempotents} 
The rectangular 0--graded object $A=k^n \in \Gamma_0$,
consisting of $n$ rows with $k$ cells in each,
plays a key role in our construction of 
1--graded idempotents. Let us define
 $P_\pm=\frac{1}{2}(1\pm A)$.

Let us enumerate  1--graded simple objects by partitions consisting of
$n$ non--increasing half--integers. The partition $S=(1/2,...,1/2)$
is used for  the object $\bullet$.\footnote{The relation between such 
partitions and dominant weights of $B_n$ is explained in the appendix.} 
The first step in the recursive construction of 
minimal 1--graded idempotents is 
given by the following proposition.

\begin{pro}\label{first}
Let $\l=(\l_1)$ be a one row Young diagram with $1\leq\l_1\leq 2k$.
The  tangles 
$$\tilde P_+(\l)=\psdiag{6}{18}{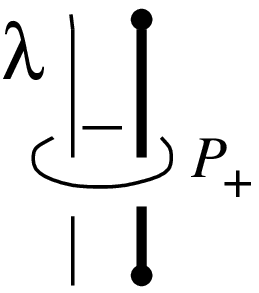} \;\;\;\;\;\;\;\;\;\;
\tilde P_-(\l)=\psdiag{6}{18}{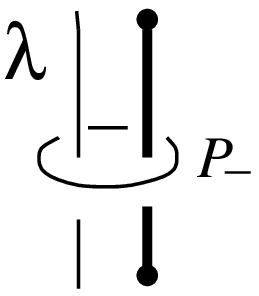}
$$
are
minimal idempotents projecting into simple objects
$(\l_1+1/2, 1/2,...,1/2)$ and $(\l_1-1/2,1/2,...,1/2)$.
Here $\tilde P_+(\l)$  projects into the first partition if $\l_1=0 \mod 2$,
otherwise into the second.

\end{pro}

\begin{proof}
Let us  
 decompose the identity of $\l\otimes S$ as follows.
\be\label{1}\psdiag{6}{18}{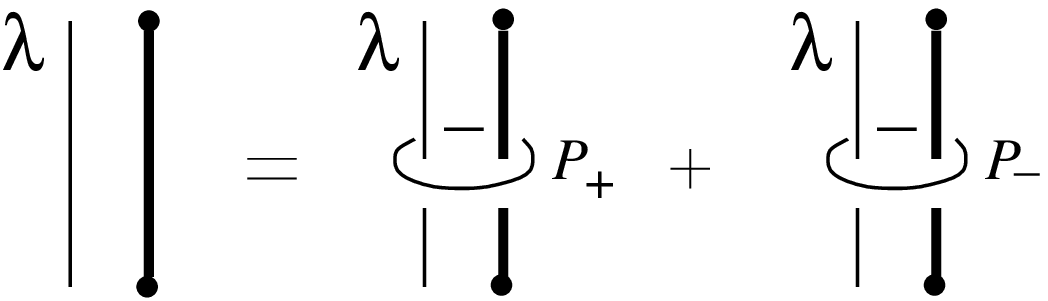}\ee
From the representation theory of $B_n$ we know that
$$ \l\otimes S= (\l_1+1/2,1/2,...,1/2)\oplus (\l_1-1/2,1/2,...,1/2)\, .$$
Therefore, dim$\,End_{{\cal T}_q (\tilde B_n)}(\l\otimes S)$ is 
maximal two. Note that $J\otimes S$ is simple, because $J=(2k+1)$ 
is invertible.
Using the equality of the colored link invariants in
${\cal T}_q(\tilde B_n)$ and $U_q(B_n)$,
semi--simplicity of the modular category for $U_q(B_n)$
and Lemma \ref{kn} in the Appendix,
we see that 
$$\tilde P_\pm(\l)\tilde P_\pm(\l)=\tilde P_\pm(\l)
\;\;\;\;\;\;\;\;\;\; \tilde P_\pm(\l)\tilde P_\mp(\l)=0$$
for any $\l\in\Gamma_0$. Furthermore, these morphisms are non--negligible.  
The claim follows.
\end{proof}

\begin{rem} In the proof we use the quantum group formulas
for the quantum dimension and the $S$--matrix.
These formulas can also be obtained
by applying an appropriate  weight system to the Kontsevich 
integral of the unknot and of the Hopf link, 
which were calculated recently in \cite{LBT}.
This will make our approach completely independent from the quantum group 
theoretical one.
\end{rem}

Let $\l=(\l_1,...,\l_p)$  be a Young diagram with $\l_1+\l_2\leq 2k$
and $2\leq p\leq n$. Let 
$\n$ be obtained by removing one cell from the last row of $\l$.
Let us assume per induction
that we can construct an idempotent $p_\m^{\n S}$
projecting  $\n\otimes S$ into a simple component $\m$.
  We know from the
classical theory that the  tensor product $\l\otimes S$
 decomposes into simple objects as follows:
\be\label{*}
 \l\otimes S=\oplus (\l_1\pm 1/2,\l_2\pm 1/2,...,\l_p\pm 1/2,1/2,...,1/2)\ee
Contributions not corresponding to non--increasing partitions
do not appear in this decomposition.
A quasi--idempotent $\tilde p^{\l S}_\m$ projecting
to a  partition $\m$ from the set 
$I=\{(\l_1\pm 1/2,\l_2\pm 1/2,...,\l_p- 1/2,1/2,...,1/2)\} $
can be obtained as follows:
$$(\tilde y_\l\otimes id_S)(id_{|\l|-1}\otimes \tilde P_+)p^{\n S}_\m
(id_{|\l|-1}\otimes \tilde P_+)(\tilde y_\l\otimes id_S)\; .$$
Here 
where $\tilde y_\l$ is the 0--graded minimal idempotent defined in \cite{BB1},
 $\tilde P_+$  is the idempotent $p^{VS}_S$
given by encircling the
2--colored line  and the 1--colored 
 line starting from the last cell in the last row   of $\l$
with a line colored by $P_+$.
Normalizing (if necessary)
this quasi--idempotent we get $p^{\l S}_\m$. 
The projector onto $(\l_1+1/2,...,\l_p+1/2, 1/2, ..., 1/2)$ 
is given by $\tilde y_\l \otimes id_S - \sum_{\m\in I} p^{\l S}_\m$.
It remains to show that $p^{\l S}_\m$ is not negligible. The trace 
of the morphism
 $(\tilde y_{(1,1)} \otimes id_S)(id_V \otimes \tilde 
P_+)$ is nonzero. Here we use that $\dim  S\neq 0$ (see next subsection). 
Analogously, the trace of
$(id_V \otimes \tilde 
P_+)(\tilde y_{(2)} \otimes id_S)$ is nonzero.
Therefore, 
$p^{\l S}_\m$ is a  composition of non--negligible  morphisms.

As a result, we can
construct 1--graded simple objects numbered by
partitions consisting of $n$ non--increasing half--integers
from the set $$\Gamma=\{ \l: \l_1+\l_2\leq 2k+1 \}.$$

We hope to be able to prove the following statement in the future.
\begin{con}
Let $\l=(\l_1,...,\l_p)$ be a Young diagram
with $\l_1+\l_2\leq 2k$, $p\leq n$.  
For  $b=(\l_1+1/2,...,\l_p+1/2, 1/2,...,1/2)$ we have
$$p^{\l S}_b=(\tilde y_\l\otimes id_S)
(\tilde P_1(\l_1)\otimes id) ...(id \otimes 
\tilde P_p(\l_p))(\tilde y_\l\otimes id_S)\, ,$$
where 
$$\tilde P_i(\l_i)=\left\{\begin{array}{rcl}
\tilde P_+(\l_i)&:& \l_i=0\mod 2\\
                  \tilde P_-(\l_i)&:& \l_i=1\mod2\\ \end{array} \right .$$

\end{con}

An example of such projection onto $(7/2, 5/2, 3/2)$
for $\l=(3,2,1)$ is drawn  below.
$$\psdiag{17}{51}{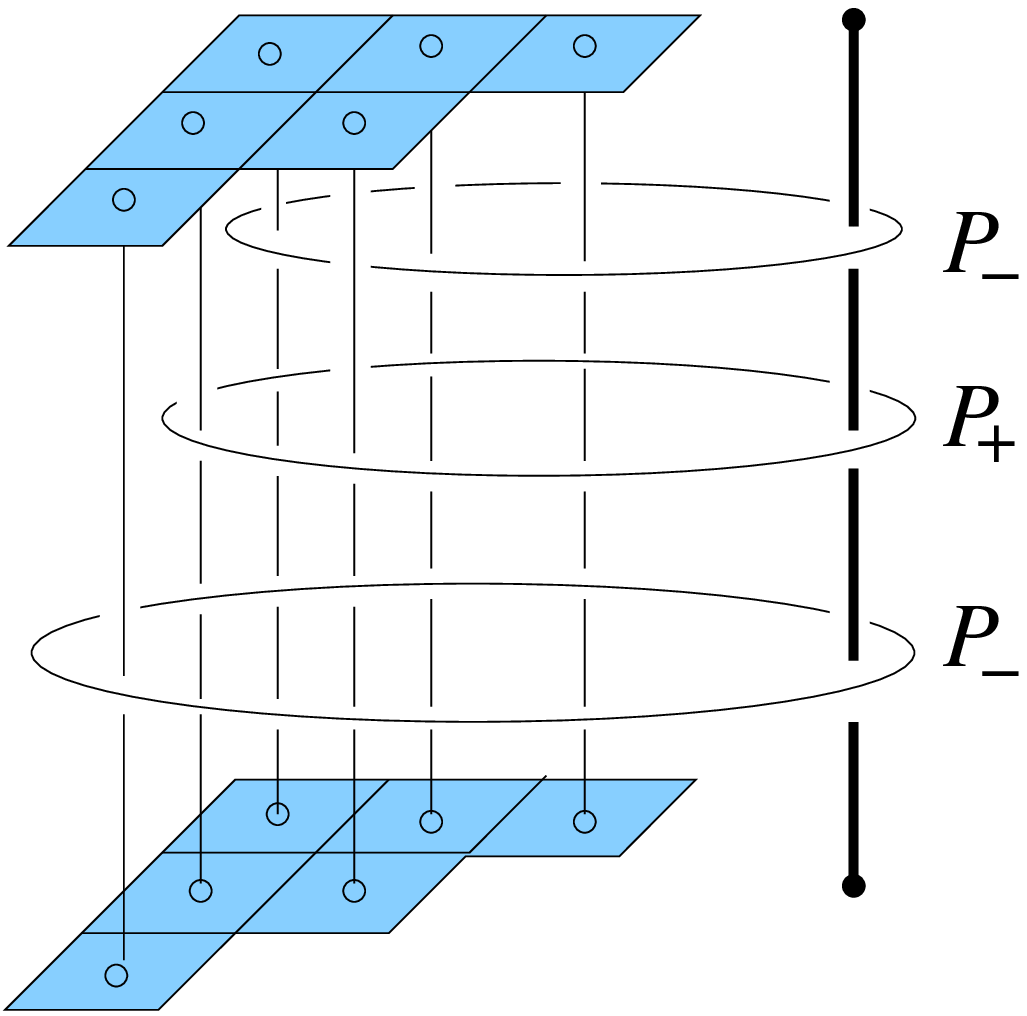}$$

\subsection{Quantum dimensions.}  
By applying the $(B_n,S)$ weight system to the Kontsevich integral of
 the $2$--colored unknot we get\footnote{This computation will 
be published elsewhere.}
$$\dim  S=(v+v^{-1})(v^3+v^{-3})...(v^{2n-1}+v^{-2n+1})\, $$
Here $q=v^2$. 
By closing (\ref{1}) with $\l=(1)$ we get
$\dim  V \;\dim  S= \dim  S +\dim  X$, which allows to calculate $\dim  X$.
We conclude that it is given by formula (\ref{dim})
in the Appendix  with $\l=(3/2,1/2,...,1/2)$. 

\begin{pro}
The quantum dimension of a simple object   $\l\in\Gamma$
is  given by formula (\ref{dim}).
\end{pro}

\proof
For integer partitions, the claim was proved in \cite{BB}.
In fact, (\ref{dim}) coincides with the formula
given in Proposition 3.3 of \cite{BB}. 
Let us assume per induction that the quantum dimensions of 1--graded
simple objects are 
given by this formula.   We are finished 
if we can show that for a $p$ row Young diagram $\l$ 
$$ \dim  S \;\dim  \l=\sum_{\tilde s\in {\mathbb Z}^p_2} \dim  (\l+\tilde s)\, ,$$
where ${\mathbb Z}^p_2=\{(\pm 1/2,...,\pm 1/2, 1/2,...,1/2)\} $.
Note that if $\l_i=\l_{i+1}$,
$$\dim  (\l_1\pm \frac{1}{2},...,\l_i-
\frac{1}{2},\l_{i}+ \frac{1}{2},...,\l_n\pm \frac{1}{2})=0\, ,$$
because terms in (\ref{dim}) corresponding to $w$ and $\sigma_i w$
cancel with each other. Here $\sigma_i$ interchanges the $i$th and $(i+1)$th
coordinates.
Using
 $$\dim  S=\sum_{\tilde s\in {\mathbb Z}^n_2} v^{2(\tilde s|\rho)}$$
we get 
\begin{align*}
 \dim  S\; \dim  \l\  =&\ \frac{1}{\psi}
\sum_{w\in W, \tilde s\in {\mathbb Z}^n_2}
sn(w) v^{2(w^{-1}(\l+\rho)+\tilde s|\rho)}  \\
=&\ \frac{1}{\psi}\sum_{w\in W, w'\in {\mathbb Z}^n_2}sn(w) v^{2(w^{-1}(\l+\rho)+ w^{-1} w'(S)|\rho)}  \\
=&\  \frac{1}{\psi}\sum_{w'\in {\mathbb Z}^n_2}\sum_{w\in W}sn(w) v^{2(\l+\rho+ w'( S)|w(\rho))}\\ =&\
\sum_{\tilde s\in {\mathbb Z}^p_2} \dim  (\l+\tilde s)\tag*{\qed} 
\end{align*}

\subsection{Twist coefficients}  
Let us denote by $t_\m$ the twist coefficient of the simple object $\m$.
Then by twisting (\ref{*}) we have the following identity
$$t_\l t_S\psdiag{6}{18}{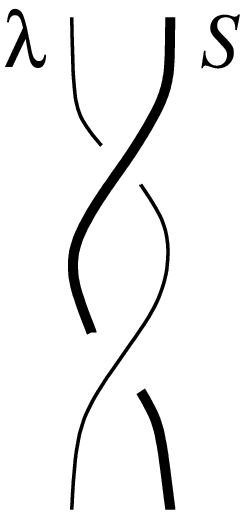}=\sum_{\tilde s\in {\mathbb Z}^p_2} 
t_{\l+\tilde s} \psdiag{6}{18}{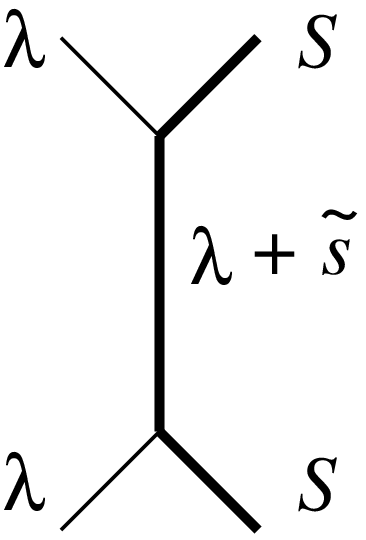}\; . $$
Replacing the positive twist with the negative one we get
a similar identity involving inverse twist coefficients.
By closing the $\l$--colored line in these two identities
we obtain
$$t_\l t_S \psdiag{6}{18}{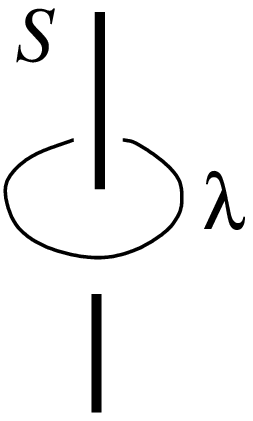}=\sum_{\tilde s}t_{\l+\tilde s}
\frac{\dim (\l+\tilde s)}{\dim \, S}\psdiag{6}{18}{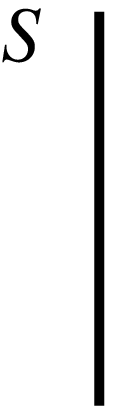}\; ,$$
$$t^{-1}_\l t^{-1}_S \psdiag{6}{18}{sl}=\sum_{\tilde s}t^{-1}_{\l+\tilde s}
\frac{\dim (\l+\tilde s)}{\dim \, S}\psdiag{6}{18}{s_bel}\; .$$
This implies the following formula:
\be\label{tw}
\sum_{\tilde s\in {\mathbb Z}^n_2} (t^{-1}_S t^{-1}_\l t_{\l+\tilde s} -
t_S t_\l t^{-1}_{\l+\tilde s})\; \dim (\l+\tilde s)=0\ee
Using (\ref{tw}) and the formulas for the quantum dimension
we can calculate twist coefficients of simple 1--graded objects recursively.
Note that $t_S=v^{n^2+n/2}$ is determined by the action of
the Casimir on the spin representation.

\begin{pro}
The twist coefficient of a simple object $\m \in \Gamma$
is given by the following formula:
\be\label{twfor}
 t_\m=v^{(\m+2\rho|\m)}\; \ee 
\end{pro}
\begin{proof}
In \cite{BB} it was shown that the twist coefficients 
of 0--graded objects are given by this formula.
Now let us assume per induction that this formula holds for 
1--graded objects. We are finished if we can prove 
(\ref{tw})  with this induction hypothesis.

Substituting (\ref{twfor}) and quantum dimensions in (\ref{tw}), we get:
\begin{eqnarray}\label{two}
v^{-(4\rho|s)} \sum_{\tilde s\in {\mathbb Z}^n_2} \sum_{w\in W}
sn(w)v^{2(\l+\tilde s+\rho|w(\rho))} v^{2(\l+\rho|\tilde s)}
\nonumber\\
 = \sum_{\tilde s\in {\mathbb Z}^n_2} \sum_{w\in W}
sn(w)v^{2(\l+\tilde s+\rho|w(\rho))} v^{-2(\l+\rho|\tilde s)}\, 
\end{eqnarray}

Using the fact that the Weyl group $W$ is a semi--direct product
of the symmetric group $S_n$ and 
${\mathbb Z}^n_2$ (acting on ${\mathbb R}^n$  by changing signs
of coordinates),
we  write $w=g'\sigma$
and $\tilde s=g(S)$  with $g, g'\in {\mathbb Z}^n_2$ and $\sigma\in S_n$.
With this notation, (\ref{two}) follows from the following two
identities:
\begin{eqnarray}\label{3}
\sum_{g,g'\in {\mathbb Z}^n_2\;, \sigma\in S_n}
sn(g'\sigma)v^{2(\l+\rho|g'\sigma(\rho)+g(S))} v^{\pm 2(g'g(S)-S|
\sigma(\rho))}\nonumber\\
=\sum_{g\in {\mathbb Z}^n_2\;, \sigma\in S_n}
sn(g\sigma)v^{2(g(\l+\rho)|\sigma(\rho)+ S)}\end{eqnarray}

\noindent
To get the second one we replace $\tilde s$ by $-\tilde s$.
In the rest of the proof we will show (\ref{3}). The idea is that 
terms with $g\neq g'$ cancel in pairs. Let us first consider the 
simplest case, when $g'(x)$ differs from $g(x)$ only 
by a sign of the $i$th coordinate. We write
$g'=g_i g$. Let us denote by $t$ the $i$th half--integer coordinate of
$\sigma(\rho)$, i.e.\ $(\sigma(\rho))_i=t$. Then there are two possibilities:
(a) there exists $j$ with $(\sigma(\rho))_j=t-1$ or (b)
$t=1/2$. In the first case, we put
$\tilde\sigma=\sigma_{ij}\sigma$ and $\tilde g=g_i g_j g$,
where $\sigma_{ij}$ interchange the $i$th and $j$th coordinates.
Analyzing the four possibilities $g(x_i)=\pm x_i$, $g(x_j)=\pm x_j$,
we see that
$$g'\sigma(\rho)+g(S)=g'\tilde\sigma(\rho)+\tilde g(S)\, .$$
The claim then follows  from the fact that
$sn(g'\sigma)=-sn(g'\tilde\sigma)$ and $(g'g(S)-S|\sigma(\rho))=
(g'\tilde g(S)-S|\tilde \sigma(\rho))$.

In case (b), we put $\tilde g'=g_i g'$, $\tilde g=g_i g$
and $\tilde \sigma=\sigma$. Case by case checking shows that
terms in (\ref{3}) corresponding to $g,g',\sigma$ and
$\tilde g,\tilde g', \tilde\sigma$ cancel with each other.
Note that
if $g'(x)$ and $g(x)$ are different for all $n$ coordinates, 
then we can proceed as in case (b). 

Let us assume that $g'(x)$ and $g(x)$ differs in less than $n$ coordinates.
Then there exists $j$ with $(\sigma(\rho))_j=1/2$. If
$g(x_j)=-g'(x_j)$, then we finish with (b),
if not, we compare $g'(x_i)$ and $g(x_i)$
with $(\sigma(\rho))_i=3/2,5/2,...$. Proceeding in this way we 
will find a pair of indices $i,j$, such that $(\sigma(\rho))_i -
(\sigma(\rho))_j=1$, $g(x_j)=g'(x_j)$, but $g(x_i)=- g'(x_i)$. 
Then we continue as in case (a).
\end{proof}

\subsection{Modular category $B^{k}_n$} 
The previous results imply that the category
${\cal T}_q (\tilde B_n)$ defined on a $(4n+4k)$th root of unity
$q$  is pre--modular. Its  simple objects
are numbered by integer or  half--integer partitions 
$\l=(\l_1, ...,\l_n)$ with $\l_1\geq...\geq\l_n\geq 0$ from
$\Gamma=\{\l:\l_1+\l_2\leq 2k+1\} .$
Let us call this category $B^{k}_n$.
\begin{thm}\label{md}
The category $B^{k}_n$ is modular.
\end{thm}
\begin{proof}
It remains to prove that $B^{k}_n$ has no 
nontrivial transparent objects.
Let $b=(b_1+1/2,b_2+1/2,..., b_n+1/2)$ be a 1--graded simple object. 
The object $J\otimes b$ is simple and is given by partition
$b'=(2k+1-b_1-1/2,b_2+1/2, ...,b_n+1/2)$. This is because
$J$ is invertible and $b'$ is the only object in $\Gamma$
with the correct twist coefficient and quantum dimension.
It follows that
$$t_J t_b\psdiag{6}{18}{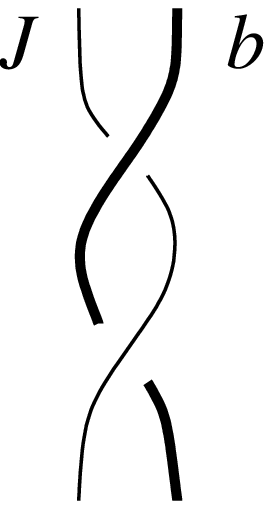}=t_{b'} \psdiag{6}{18}{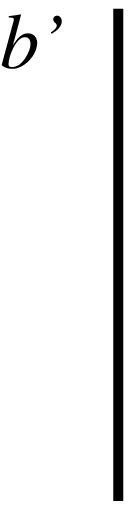}\; . $$
Inserting  twist coefficients  we obtain that
the braiding coefficient of $J$ and $b$ is $(-1)$. This implies that $J$
is not transparent in $B^{k}_n$ and that no 1--graded simple object
can be transparent. But the 0--graded part of $B^{k}_n$
has not even a further nontrivial $0$--transparent object.
\end{proof}
\subsection{Refinements}
It was shown by Blanchet in \cite{B} that any 
modular category with an invertible object $J$ of order 2 (i.e.\ $J^2=1$),
whose twist coefficient is $(-1)$ and quantum dimension is $1$,
provides invariants of 3--manifolds equipped with spin structure.
The  ${\mathbb Z}_2$--grading defined in \cite{B} 
by means of $J$  coincides with the  one used in this paper.
 The Kirby color
decomposes as $\Omega=\Omega_0+\Omega_1$ 
according to this grading.
The  invariants 
of closed 3--manifolds equipped with spin structure
are defined by putting the 1--graded 
Kirby color on the components of a surgery link belonging to
the so--called  characteristic sublink (defined by the
 spin structure) and the 0--graded Kirby color on the other components.
The ordinary 3--manifold invariant decomposes into a sum of refined invariants 
over all spin structures. A spin TQFT can also be constructed (see e.g.\
\cite{b}).
It associates a vector space $V(\Sigma_g,\sf s)$ to a genus $g$
surface $\Sigma_g$ with spin structure $\sf s$.

\begin{pro}
The category $B^{k}_n$ provides a spin TQFT.
Furthermore,
\begin{eqnarray*}
 \frac{4^{g}}{\la \Omega \ra^{g-1}}\;\;
\dim V(\Sigma_g, {\sf s})& = & 
 \sum_{\l\in \Gamma \setminus \Gamma_1} (dim\; \l)^{2-2g}\\
& & +\;\;\; (-1)^{{\rm Arf}({\sf s})}
\;\;{2^{g}}\;\;
 \sum_{\l\in \Gamma_1} (dim\; \l)^{2-2g}\end{eqnarray*}
where $\la \Omega \ra$ is the invariant of 
the Kirby--colored unknot, ${\rm Arf}(\sf s)$ is the Arf invariant
and $\Gamma_1=\{ \l\in\Gamma: \l_1=k+1/2\}$. 
\end{pro}
\begin{proof}
A spin Verlinde formula was computed by Blanchet
in \cite[Theorem 3.3]{B}. It uses
 the action of $J$ on $\Gamma$ given by the tensor product. 
In our case,
 $$J\otimes (\l_1,\l_2,...,\l_n) =(2k+1-\l_1,\l_2,...,\l_n)$$
(compare  \cite{BB} and the proof of Theorem \ref{md}.) 

Therefore, there are only two different  cases.
If $\l_1\neq k+ 1/2$, then 
$\#orb(\l)=2$, $|Stab(\l)|=1$. If $\l_1= k+ 1/2$, then
$\#orb(\l)=1$, $|Stab(\l)|=2$.
The result follows by the direct application of the Blanchet formula.
\end{proof}

\section{Even orthogonal modular categories}
In this section we construct even orthogonal modular categories.
We show that corresponding invariants admit  cohomological 
refinements  and calculate the refined Verlinde formulas.

\subsection{0--graded objects}
Let us use the standard notation $D_n$
for ${\mathfrak so}_{2n}$. We fix a primitive $(2k+2n-2)$th
root of unity $q$  and $v$ with $v^2=q$. According to \cite{BB},
the category ${\cal T}_q (D_n)$
has  the following set of 0--graded simple objects 
$$\tilde
\Gamma=\{\l:\l_1+\l_2\leq 2k, \l^\vee_1+\l_2^\vee\leq 2n\}\; .$$  

This set contains two invertible
objects of order two: $1^{2n}$ and $2k$.
They are  $0$--transparent, with
 twist coefficients and quantum dimensions are equal to $1$. 
This implies that the 0--graded part of ${\cal T}_q (D_n)$ is modularizable.
After modular extension by the group generated by  $1^{2n}$
we get  a new category, which will be denoted by  ${\cal T}_q (\tilde D_n)$.
  The objects $  1^{2n}\otimes \n$ and $\n$
are isomorphic there for any
$\n\in \tilde \Gamma$. The 0--graded objects $\l$
with $\l^\vee_1=n$ do not remain  simple  in  ${\cal T}_q (\tilde D_n)$
and  decompose as  $\l=\l_-+\l_+$.
The objects $\l_\pm$ have the same quantum dimensions
 and twist coefficients.
We will use
the partitions $(\l_1,...,\pm\l_n)$ for $\l_\pm$. 
The set of 0--graded simple objects of  ${\cal T}_q (\tilde D_n)$ is
$$\Gamma_0=\{\l:\l_1+\l_2\leq 2k, \l^\vee_1<n\}\; \cup
\{\l_\pm:\l_1+\l_2\leq 2k, \l^\vee_1=n\}\; .$$

Let  $A_+$  be the  simple 0--graded object of  ${\cal T}_q (\tilde D_n)$ 
obtained after splitting 
 of $A=(k,k,...k)=k^n$.
Let $i=v^{n+k-1}$ and
   $P_\pm=\frac{1}{2}(1\pm (-i)^n A_+)$. 
Then analogously to the odd orthogonal case, encircling of a spinor 
$b=(\frac{2b_1+1}{2},...,\frac{2b_n+1}{2})$ by $P_+$ gives 
the identity morphism if $\sum_i b_i=0\mod 2$ and is zero 
otherwise. 
The proof is given in the appendix.

\subsection{Recursive construction of 1--graded idempotents}
The group  $D_n$ has two spin representations
$S_\pm$ given by the highest weights $(1/2,...,\pm1/2)$.
In order to distinguish them we put an orientation on the  2--colored lines.

Let $w\in {\mathbb Z}_2^{p}$ act on coordinates of  ${\mathbb R}^{p}$ by 
sign changing. We put $sn(w)=1$ if it changes the signs of an even number 
of coordinates and $sn(w)=-1$ otherwise.
Let $s=(1/2,...,1/2)\in {\mathbb R}^{p}$.
For any highest weight 
 $\l=(\l_1,...,\l_p,0,...,0)$, 
the tensor product $\l\otimes S_\pm$ decomposes in 
$D_n$ as follows:
$$\l\otimes S_\pm=
\bigoplus_{w\in {\mathbb Z}^p_2}
(\l_1+w(s_1),..., \l_p+w(s_p),1/2,..., \pm sn(w)1/2)$$
If $p=n$, we have
$$\l\otimes S_\pm=
\bigoplus_{w\in {\mathbb Z}^{n-1}_2}
(\l_1+w(s_1), \l_2+w(s_2),...,\l_n \pm sn(w)1/2)\, .$$
In particular,
$$V\otimes S_+=S_- + (3/2,1/2,...,1/2)\; ,\;\;\;\;
V\otimes S_-=S_+ + (3/2,1/2,...,-1/2)\; .$$
The corresponding idempotents are given by $\tilde P_\pm(V)$.

Suppose that we can decompose $\n\otimes S_\pm$ 
into simple objects
if $\n$ is obtained
by removing one cell from the last row of $\l$. 
Then the projection $p^{\l S_+}_\m$
to $\m\subset \l\otimes S_+$  (with 
$\m\neq  h=(\l_1+1/2,...,\l_k+1/2,1/2,...,1/2)$)
can be obtained by normalizing the  following morphism
$$(\tilde y_\l\otimes id_{S_+})(id_{|\l|-1}\otimes \tilde P_+)p^{\n S_-}_\m
(id_{|\l|-1}\otimes \tilde P_+)(\tilde y_\l\otimes id_{S_+})\; ,$$
where  $\tilde P_+$  is given  by encircling the
2--colored line  and the 1--colored 
 line starting from the last cell in the last row  of $\l$
with a line colored by $P_+$.
The idempotent  $p^{\l S_+}_h$ is given by
$\tilde y_\l\otimes id_{S_+} -\sum_\m p^{\l S_+}_\m$. The case $\l\otimes S_-$ 
is similar.

\subsection{Modular category $D^{k}_n$}
Analogously  to the  odd orthogonal case, 
one can  show that the
quantum dimensions of spinors are given by the formula (\ref{dim})
and
the twist coefficient of a simple object $\m$ of ${\cal T}_q (\tilde D_n)$
is $v^{(\m+2\rho|\m)}$.

We conclude that the category  ${\cal T}_q (\tilde D_n)$ at 
a $(2n+2k-2)$th root of unity $q$ 
is pre--modular.   Its  simple objects are given by 
integer or half--integer
partitions $\l=(\l_1,...,\l_{n-1},\pm \l_n)$ with
$\l_1\geq\l_2\geq ...\geq \l_n\geq 0$ from the set
$\Gamma=\{\l:\l_1+\l_2\leq 2k\}$.
Let us denote this  category by $D^{k}_n$.
Taking into account that the object $2k$ has the braiding coefficient $(-1)$
with any spinor, we derive that $D^{k}_n$ is modular.

\subsection{Refinements}
The category $D^{k}_n$ has an invertible  object $J$ of order 2,
whose twist coefficient and quantum dimension are equal to 1.
It was shown in \cite{B} that any such
modular category  
provides an invariant of a 3--manifold $M$  equipped with
a first ${\mathbb Z}_2$--cohomology class. 
More precisely, the object $J$  defines a grading in the category,
which coincides with the ${\mathbb Z}_2$--grading used above.
For a closed 3--manifold $M$,
any $h\in H^1(M,{\mathbb Z}_2)$ can be represented by a sublink
of a surgery link for $M$ belonging to the kernel 
of the linking matrix modulo $2$. The invariant of a pair $(M,h)$
is then defined by putting 1--graded Kirby colors on this sublink
and 0--graded ones on the other components. This
construction can be extended
to manifolds with boundary (see \cite{b})
and  leads to a cohomological  TQFT (compare \cite{LT}). 

\begin{pro}
The category $D^k_n$ leads to a cohomological TQFT.
The dimension of the TQFT module associated with  a pair
$(\Sigma_g,h)$, $h\in H^1(\Sigma_g)$,  is given by the following 
formulas. For  $h\neq 0$,
\begin{eqnarray*}
\dim V(\Sigma_g, h) =  
 \frac{\la \Omega \ra^{g-1}}{4^{g}}\;\; \sum_{\l\in \Gamma \setminus \Gamma_1} (dim\; \l)^{2-2g}\\
\end{eqnarray*}
where $\Gamma_1=\{\l\in \Gamma: \l_1=k, \l^\vee_1<n  \}$. For $h=0$,
$$
\dim V(\Sigma_g, h) =  \frac{\la \Omega \ra^{g-1}}{4^{g}}\;\;
\left( \sum_{\l\in \Gamma \setminus \Gamma_1} (dim\; \l)^{2-2g}
+ 4^g \;\;\sum_{\l\in \Gamma_1} (dim\; \l)^{2-2g}\right)
$$

\end{pro}
\begin{proof}
In  \cite[Theorem 5.1]{B} Blanchet give a refined Verlinde 
formula for cohomological TQFT's.
In our case,
 $J$ acts on  $ \Gamma$ as follows:
 $$J\otimes (\l_1,\l_2,...,\l_n) =(2k-\l_1,\l_2,...,-\l_n)$$
This is because,
$\l\otimes J$ is simple and 
it contains the  object from the right hand side of the above formula
by classical representation theory.

Therefore, we have only two cases. If $\l\in \Gamma\setminus
\Gamma_1$, then 
$\#orb(\l)=2$, $|Stab(\l)|=1$. If $\l\in \Gamma_1$, then
$\#orb(\l)=1$, $|Stab(\l)|=2$.
The result follows by the direct application of the Blanchet formula.
\end{proof}

\section{Relation with quantum groups}
 We show that our orthogonal modular categories
  are equivalent
to the quantum group theoretical ones. Further, we compare
results about refinements and  level--rank duality.

\subsection{Equivalence}
Let us call two modular categories {\it equivalent} if there exists a 
bijection between their sets of simple objects providing an equality 
of the corresponding colored link invariants.
This implies that the associated TQFT's are isomorphic (see 
\cite[III, 3.3]{T}).

\begin{thm}
{\rm i)}\qua The category $B^k_n$ is equivalent to the
modular category defined for $U_q(B_n)$  at a $(4n+4k)$th root of unity $q$.

{\rm ii)}\qua  The category $D^{k}_n$ is equivalent to the
modular category defined for $U_q(D_n)$  
at a $(2n+2k-2)$th root of unity $q$.
\end{thm}

\begin{proof}
The construction of modular categories from quantum groups is given in
\cite{Le}. Let us recall the main results. 

Let  $\mathfrak g$ be  a finite--dimensional simple
 Lie algebra over $\mathbb C$.
Let $d$ be the maximal absolute value of the non--diagonal 
entries of its Cartan matrix.
Let us denote by $C$ the set of the dominant weights of 
$\mathfrak g$. We normalize the 
inner product $(\;|\;)$  on the weight space, such that  the square 
length of any short root is 2. Finally, we denote by $\beta_0$ 
 the long root in  $C$.

The quantum group $U_q(\mathfrak g)$ at a primitive root of unity
$q$ of order $r$ provides a modular category if $r\geq dh^\vee$,
where $h^\vee$ is the Coxeter number.  
 This modular category has the following
set of simple objects.
$$C_L=\{x\in C: (x|\beta_0)\leq d L\}$$
Here $L:=r/d-h^\vee$ is the level of the category.
 
{\rm i)}\qua Let $\mathfrak g=B_n$. Then we have 
 $\beta_0=(1,1,0,...,0)$ in the basis
chosen in the appendix, $d=2$ and  $h^\vee=2n-1$. For $r=4n+4k$,
 $C_L$ is in bijection  with the set 
$\Gamma=\{\l:\l_1+\l_2\leq 2k+1\}$ of simple objects of $B^k_n$.
Moreover, this bijection induces an equality
of the colored link invariants 
due to the result of Le--Murakami  \cite{LM}.

{\rm iI)}\qua Let $\mathfrak g=D_n$. Then  
 $\beta_0=(1,1,0,...,0)$ in the basis
chosen in the appendix, $d=1$ and  $h^\vee=2n-2$.
For $r=2n+2k-2$,
  $C_L$ coincides  with the set 
 of simple objects of $D^k_n$. The claim follows then 
as above from  \cite{LM}.
\end{proof}

\subsection{Refinements}
Cohomological refinements 
in categories obtained from quantum groups were
studied in \cite{LT}. 
For  type $D$, Le and Turaev consider
 cohomology classes  with
coefficients in ${\mathbb Z}_4$ or ${\mathbb Z}_2\times {\mathbb Z}_2$.
The statements about existence of spin refinements
in modular categories of type $B$  and
about ${\mathbb Z}_2$--cohomological refinements for type $D$ 
seem to be new. 

\subsection{Level--rank duality}

It was shown in \cite{BB} that the categories $B^{n,-k}$
and $D^{n,k}$ have their
level--rank dual partners. For
quantum groups this means the following.
Let us denote by $\tilde D^{n,k}$
 the modular category for
$U_q(D_n)$ at a $(2n+2k-2)$th root of unity quotiented
by spinors and the action of the transparent object $2k$.
Then $\tilde D^{n,k}$ 
is isomorphic to $\tilde D^{k,n}$.
  The isomorphism is given
by sending $v$ to $-v^{-1}$  and by `transposing' the partitions.
(It is helpful to use the geometric approach to modularization
functors in order to construct the isomorphism.)
In the odd orthogonal case this duality does not exist
on the quantum group level, because the corresponding quotients can 
not be constructed. Transparent objects have twist coefficients $(-1)$.

Unfortunately, this
level--rank  duality between 0--graded parts of
$B^{k}_n$ and $D^{k}_n$ does not extend to the full
categories.
Even the cardinalities of the sets of simple objects in $B^{k}_n$
and $B^{n}_k$ as well as in $D^{k}_n$ and $D^{n}_k$
 are different in general. This suggests the existence of 
bigger   categories 
with more symmetric sets of objects admitting 
level--rank duality. One possibility to construct them would be
to take the 0--graded part of $B^k_n$
and to add two different spin representations
using the $B_n$ and $B_k$  weight systems. This will be
studied in the forthcoming  paper with C. Blanchet.

\section{Appendix}
\subsection{Odd orthogonal case}
Let $\{e_i\}_{i=1,2,...,n}$ be the standard base of ${\mathbb R}^n$ 
with the scalar product $(e_i|e_j)=2\delta_{ij}$. Any weight of $B_n$ has
all integer or all half--integer coordinates in this base. 
We write $\l=(\l_1,...,\l_n)$ if $\l=\sum^n_{i=1} \l_i e_i$.
Any weight $\l$ with
 $\l_1\geq\l_2\geq...\geq \l_n\geq 0$ is a highest weight
of an irreducible representation of $B_n$ or a dominant weight.
 The half sum of all positive roots of $B_n$ we denote by
$\rho=(n-1/2, n-3/2,...,1/2)$. 

Let us consider the quantum group $U_q(B_n)$, where $q=v^2$ 
is a primitive  $(4n+4k)$th root of unity. The set 
of simple objects (or dominant weights)
of the corresponding modular category is
  $\Gamma=\{\l: \l_1+\l_2\leq 2k+1 \}$, where $\l$  is a highest weight
of $B_n$ \cite{Le}.
The quantum dimension of a simple object 
 $\l\in\Gamma$ is given by the following formula:
\be\label{dim}
dim\; \l=\frac{1}{\psi}\sum_{w\in W} sn(w)v^{2(\l+\rho|w(\rho))}\ee
$$\psi=\sum_{w\in W} sn(w)v^{2(\rho|w(\rho))}=\prod_{positive\; roots\; \a}
v^{(\rho|\a)}-v^{-(\rho|\a)}\; $$
Here $W={\mathbb Z}^n_2\propto S_n$ is  the Weyl group of $B_n$
 generated by reflections on the hyperplanes orthogonal to
the roots $\pm e_i\pm e_j, \pm e_i$.
Furthermore, 
for  the invariant of the Hopf link, whose components are colored 
by   $\m,\n\in \Gamma$,
we have
\be\label{s}
S_{\m\n}=\frac{1}{\psi}\sum_{w\in W} sn(w) v^{2(\m+\rho|w(\n+\rho))}\, .\ee

\begin{lem}\label{kn}
Let  $b=(\frac{2b_1+1}{2}, \frac{2b_2+1}{2},...,\frac{2b_n+1}{2})$
be a  dominant weight with half--integer coordinates
  and $A=(k,k,...,k)$,
 then
$$\psdiag{6}{18}{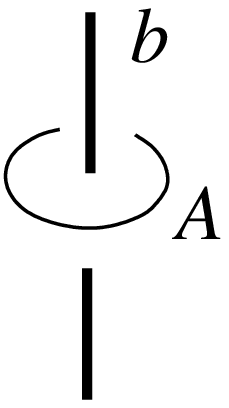}=(-1)^{b_1+b_2+...+b_n}\psdiag{6}{18}{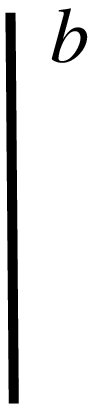} .$$ 
\end{lem}
\begin{proof}
The coefficient to determine is equal to $S_{bA}(\dim \;b)^{-1}$.
From (\ref{dim}) and (\ref{s}) we have
$$dim\; b=\frac{1}{\psi}\sum_{w\in W} sn(w)v^{2(b+\rho|w(\rho))}$$
$$S_{bA}=\frac{1}{\psi}\sum_{w\in W} sn(w) v^{2(b+\rho|w(A+\rho))}$$
Let us write
$A+\rho=(k+n)(1,1,...,1)+w_1(\rho)$, where $w_1\in W$ and $sn(w_1)=(-1)^
{n(n+1)/2}$. Then from  $v^{4n+4k}=-1$ we have
$$v^{2(k+n)\sum_{i}(b+\rho|\pm e_i)}
=(-1)^{(b_1+b_2+...+b_n+n(n+1)/2)}\, $$
or
$$sn(w)v^{2(b+\rho|w(A+\rho))}=(-1)^{b_1+b_2+...+b_n}sn(ww_1)v^{2(b+\rho|ww_1(\rho))}\, .$$
The result follows.
\end{proof}

\subsection{Even orthogonal case}
Let $\{e_i\}_{i=1,2,...,n}$ be the standard base of ${\mathbb R}^n$ 
with the scalar product $(e_i|e_j)=\delta_{ij}$. Any weight of $D_n$ has
all integer or all half--integer coordinates in this base. 
Any weight $\l=(\l_1,...,\pm \l_n)$  with
 $\l_1\geq\l_2...\geq \l_n\geq 0$ is a highest weight
of an irreducible representation of $D_n$. 
 The half sum of all positive roots of $D_n$ we denote by
$\rho=(n-1, n-2,...,1,0)$.

Let us consider the quantum group $U_q(D_n)$, where 
 $q$ is a primitive  $(2n+2k-2)$th root of unity and $v^2=q$. 
The set of simple objects (or dominant weights)
 of the corresponding modular category is
 $\Gamma=\{\l: \l_1+\l_2\leq 2k  \}$. 
The formulas (\ref{dim}) and (\ref{s}) hold for the highest weights
from $\Gamma$, where 
 the Weyl group $W$ of $D_n$ is generated by
reflections on the hyperplanes orthogonal to the roots $\pm e_i\pm e_j$.
This group contains $S_n$. The kernel of the projection to $S_n$
consists of transformations acting by $(-1)$ on an even number of axes. 

\begin{lem}
Let  $b=(\frac{2b_1+1}{2}, \frac{2b_2+1}{2},...,\frac{2b_n+1}{2})$ be
a  do\-mi\-nant weight with half--integer coordinates,
 $A_+=(k,k,...,k)$ and $i=v^{n+k-1}$,
 then
$$\psdiag{6}{18}{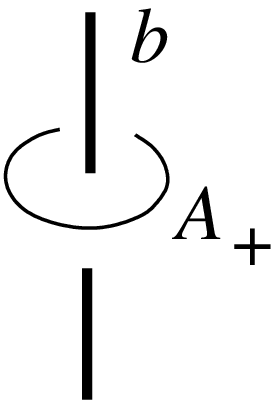}=(-1)^{b_1+b_2+...+b_n} \; i^n
\psdiag{6}{18}{b_bel} .$$ 

\end{lem}
\begin{proof}
As before,
the coefficient to determine is equal to $S_{bA_+}(\dim \;b)^{-1}$.
We have
$$S_{bA_+}=\frac{1}{\psi}\sum_{w\in W} sn(w) v^{2(b+\rho|w(A_+ +\rho))} \, .$$
Let us write
$A_+ +\rho=(k+n-1)(1,1,...,1)+w_1(\rho)$, where $w_1\in W$ and $sn(w_1)=(-1)^
{n(n-1)/2}$. Then
$$v^{2(k+n-1)(b+\rho|w(1,...,1))}=(-1)^{b_1+b_2+...+b_n+n(n-1)/2}\; i^n$$
for any $w\in W$. The result follows as in the odd orthogonal case.
\end{proof}

Note that for $A_-=(k,k,...,-k)$ an analogous statement holds:
$$\psdiag{6}{18}{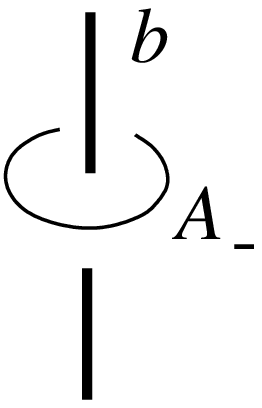}=(-1)^{b_1+b_2+...+b_n} (-i)^n
\psdiag{6}{18}{b_bel} $$

\bibliographystyle{amsplain}

\Addresses\recd
 
\end{document}